\documentclass[12pt,a4paper]{article}
\usepackage{amsmath,amsthm,pstricks,indent}

\makeatletter

\renewcommand{\section}{\@startsection
	{section}%
	{1}%
	{0mm}%
	{-\baselineskip}%
	{0.5\baselineskip}%
	{\bfseries\normalsize\centering}}%
\makeatother

\newcommand{\osmall}{\mathop{o}\nolimits}

\newcommand{\divides}{\mathrel{\vert}}

\newtheorem{teo}{Theorem}
\newtheorem{cor}{Corollary}
\newtheorem{dfn}{Definition}

\title{Enumeration of chord diagrams}
\author{A. Khruzin}
\date{}

\begin{document}

\maketitle
\makeatletter
\renewcommand{\@makefntext}[1]{\parindent=1em\noindent #1}
\makeatother
\setlength{\footnotesep}{0.8\baselineskip}
\renewcommand{\thefootnote}{}
\makeatletter
\renewcommand{\@makefntext}[1]{\parindent=1em\noindent%
              \hbox to 1.8em{\hss$^{\@thefnmark}$}#1}
\makeatother
\renewcommand{\thefootnote}{\arabic{footnote}}
\setcounter{footnote}{0}

\begin{abstract}
We determine the number of nonequivalent chord diagrams of order $n$
under the action of two groups, $C_{2n}$, a cyclic group of 
order $2n$, and $D_{2n}$, a dihedral group of order $4n$. 
Asymptotic formulas are also established.
\end{abstract}

\section{}

Given $2n$ different points on a circle ask a question: in how many different 
ways may the points be joined by chords. The answer depends of course on our
understanding of the word "different".

The configuration (actually a graph) consisting of the circle and $n$ chords
joining $2n$ different point is called a chord diagram of order $n$ or, shortly,
$n$-diagram. In the present paper, we let a group $G$ act on the circle and 
consider two $n$-diagrams as indistinguishable or equivalent if the one is 
transformed into the other by a suitable element of the group.

With the identity group acting on the circle, all $n$-diagrams are distinct
and there are altogether $\frac{(2n)!}{2^n n!} = (2n - 1)!!$ diagrams with $n$ chords.
This case is well-studied: Errera \cite{Er} (see also A.~M.~Jaglom and I.~M.~Jaglom
\cite{JJ}) determined the number of $n$-diagrams with the additional 
requirement that no chords intersect inside the circle which equals 
$\frac{(2n)!}{n!(n + 1)!}$. Touchard \cite{To} and later 
Riordan \cite{Ri} extended that result to enumeration of $n$-diagrams 
by the number of crossings of the chords, which is given by a generating 
function $T_n(x)$ satisfying the following relation
\[
(1 - x)^n T_n(x) = \sum_{j = 0}^n (-1)^j t_{nj}x^J
\]
with
\[
J = \binom{j + 1}{2}, \quad t_{nj} = \frac{2j + 1}{2n + 1} \binom{2n + 1}{n - j}.
\]

Generally chord diagrams are not strict graphs as they may have double edges.
Chord diagrams which are strict graphs were considered by Hazewinkel and 
Kalashnikov \cite{HK}. Let $b_{2n}$ be the number of strict $n$-diagrams 
under the identity group. They proved that 
$a_{2n} = \sum_{i = 1}^n b_{2i}$ satisfies the recurrence 
$a_{2n} = (2n - 1) \, a_{2n - 2} + a_{2n - 4}$.

Chord diagrams considered to within an equivalence induced
by a cyclic group action appear in different contexts. In the theory 
of knots they are used to describe Vassiliev knot invariants
in a purely combinatorial way by considering the algebra of
functions defined on the set of chord diagrams satisfying certain
linear equations \cite{CDL}, \cite{BN}. Chord diagrams also appear 
in the classification of vector fields and smooth functions on surfaces 
up to a homeomorphism where stable separatrices of critical points
together with a component of the surface boundary form a chord
diagram embedded into the surface \cite{Sh}.

\section{}

\begin{dfn} A chord diagram of order $n$ is a 3-regular graph with the 
vertex set $[2n] = \{1, 2, \dots, 2n \}$ containing the $2n$-circuit 
$\Delta_{2n} = (1 \, 2 \, \dots 2n)$ as a subgraph. The circuit is 
called a circle, the edges not belonging to the circuit are called chords.
\end{dfn}

\begin{dfn} Let a group $G$ act on the circuit $\Delta_{2n}$. 
Two $n$-diagrams $\Gamma_1$ and $\Gamma_2$ are said to be equivalent if there 
is a $g \in G$ which takes the chords of $\Gamma_1$ into the chords of $\Gamma_2$.
\end{dfn}

Our first result is the following
\begin{teo}
The number of nonequivalent $n$-diagrams under the action of a group $G$ equals
\begin{equation}\label{eq1}
\frac{1}{|S_n \wr S_2| \cdot |G|} \sum_{\pi \in S_n \wr S_2} \sum_{\eta \in G} \;
\prod_i i^{\pi_i} \: \eta_i \: (\eta_i - 1) \: \ldots \: (\eta_i - \pi_i + 1).
\end{equation}
Here $S_n \wr S_2$ is the wreath product of two symmetric groups $S_n$ and $S_2$, 
$\pi \in S_n \wr S_2$ has cycle type $1^{\pi_1} 2^{\pi_2} \ldots (2n)^{\pi_{2n}}$,
$\eta \in G$ cycle type $1^{\eta_1} 2^{\eta_2} \ldots (2n)^{\eta_{2n}}$ and
the product is taken over all $i \in [2n]$ such that $\pi_i > 0$, the product 
being equal to zero if $\pi_i > \eta_i$ for some $i$.
\end{teo}

\begin{proof} 
For each chord diagram, a subgraph consisting of its chords is a 1-factor.
The chords of all $n$-diagrams constitute the complete graph $K_{2n}$, the chords
of each single $n$-diagram being again a 1-factor of $K_{2n}$. The action of $G$
on $\Delta_{2n}$ induces an action on the set $\mathcal F_1$ of all 1-factors of 
$K_{2n}$. The orbits of $\mathcal F_1$ under that action are in a one-one 
correspondence with the nonequivalent $n$-diagrams. 

A 1-factor of $K_{2n}$ can be represented by a $n \times 2$ matrix
whose entries belong to $[2n]$ and are all distinct: each row corresponds
to an edge, two row entries being the end points of the edge. This 
correspondence is not unique. It is defined to whithin an equivalence 
induced on the set of such matrices by independently permuting entries 
in each row and permuting the rows bodily. This amounts to the action 
of the wreath product $S_n \wr S_2$. The action of $G$ on the set 
of 1-factors is equivalent to the action of $G$ on the set of 
$2 \times n$-matrices. 

We thus arrive at the following setting: given two sets $[n] \times [2]$ and 
$[2n]$, consider the set of bijective mappings $[n] \times [2] \to [2n]$.
The wreath product $S_n \wr S_2$ acts on the set $[n] \times [2]$ by the rule 
$(\tau, \bar\sigma) \cdot (i, j) = (\tau(i), \sigma_i(j))$ where
$\tau \in S_n$, $\bar\sigma = (\sigma_1, \ldots \sigma_n) \in S_2^n$, 
and the group $G$ acts on the set $[2n]$.
Two mappings $f_1, f_2 : [n] \times [2] \to [2n]$ are equivalent if there 
exist a $\pi \in S_n \wr S_2$ and an $\eta \in G$ such that
\[
f_1(\pi (i, j)) = \eta f_2 ((i, j))
\]
for all $(i, j) \in [n] \times [2]$. The equivalence classes of mappings are in 
a one-one correspondence with the orbits of $\mathcal F_1$. In this setting, an 
argument of de Bruijn \cite{dB} applies which he used to prove a theorem on 
the number of classes of bijective mappings. The proofe is complete.
\end{proof}

\section{}

We now specialize $G$ to a cyclic group $C_{2n}$ of order $2n$ and obtain 
a much simpler expression for the number of nonequivalent $n$-diagrams.
\begin{teo}\label{t2}
The number $c_n$ of nonequivalent $n$-diagrams under the action 
of a cyclic group $C_{2n}$ equals
\begin{equation}\label{eq2}
c_n = \frac{1}{2n} \sum_{i \divides 2n} \varphi(i) \nu_n(i),
\end{equation}
where $\varphi(i)$ is the Euler function and
\begin{equation}\label{eq3}
\nu_n(i) = \begin{cases}
i^{n/i} \, \left( 2n/i - 1 \right)!! & \text{$i$ odd}, \\
{\displaystyle \sum_{k = 0}^{\lfloor \frac{n}{i} \rfloor} 
\binom{2n/i}{2k} \, i^k \, (2k - 1)!!}  & \text{$i$ even}
\end{cases}
\end{equation}
for $i \divides 2n$.
\end{teo}

\begin{proof} Each permutation $\eta \in C_{2n}$ has cycle type $i^{2n/i}$, 
$i \divides 2n$. If $\pi \in S_n \wr S_2$ 
has a cycle of length $\ne i$ the product in~(\ref{eq1}) equals zero,
otherwise it reduces to a single term
$i^{\pi_i} \eta_i (\eta_i - 1) \ldots (\eta_i - \pi_i + 1) = i^{2n/i} (2n/i)!$
since $\pi_i = \eta_i = 2n/i$. The double sum in~(\ref{eq1}) can then be replaced with 
a single sum over all $i \divides 2n$. The group $C_{2n}$ contains
$\varphi(i)$ permutations of cycle type $i^{2n/i}$. Denoting by $\psi_n(i)$ the
number of permutations of the same cycle type $i^{2n/i}$ in $S_n \wr S_2$ 
we rewrite~(\ref{eq1}) as follows
\begin{equation}\label{eq4}
c_n = \frac{1}{2^n \: n! \: 2n} \sum_{i \divides 2n} 
i^{2n/i} \, (2n/i)! \, \varphi(i) \, \psi_n(i).
\end{equation}

To determine $\psi_n(i)$ we first establish a relationship between the cycles 
of $\tau \in S_n$ and the cycles of $(\tau, \bar\sigma) \in S_n \wr S_2$.

If $K_\pi \subset [n] \times [2]$ is a cycle of $\pi = (\tau, \bar\sigma)$ 
then its projection onto $[n]$ is a cycle of $\tau$.

If $K_\tau \subset [n]$ is a cycle of $\tau \in S_n$ then, for any 
$\bar\sigma \in S_2^n$, the length of a cycle of $(\tau, \bar\sigma) \in S_n \wr S_2$
induced by $K_\tau$ depends only on $\sigma_k \in S_2$, $k \in K_\tau$.
Let $L = \{ k \in K_\tau : \sigma_k \ne e\}$ where $e$ is the identity
permutation. If $|L|$ is odd then $K_\tau$ induces one cycle of length 
$2 \, |K_\tau|$. If $|L|$ is even then $K_\tau$ induces two cycles both
of length $|K_\tau|$.

Let $i \divides 2n$ be odd. For $\pi = (\tau, \bar\sigma)$ to
have cycle type $i^{2n/i}$ the permutation $\tau$ must have cycle type $i^{n/i}$.
The number of such $\tau \in S_n$ equals
\[
\frac{n!}{i^{n/i} \: (n/i)!}.
\]
Now we fix $\tau$ and count the number of $\bar\sigma \in S_2^n$ such that 
$(\tau, \bar\sigma)$ has cycle type $i^{2n/i}$. Each cycle $K_\tau$ of $\tau$ 
induces two cycles of $(\tau, \bar\sigma)$ of the same length $i$. Hence $|L|$ 
must be even for each $K_\tau$. There are
\[
\sum_{m \leq i, \; m \: \text{even}} \binom{i}{m} = 2^{i - 1}
\]
choices for $L \subset K_\tau$. Clearly, $L$ uniquely determines $\sigma_k$
for all $k\in K_\tau$. For different cycles of $\tau$ the choices
are independent, so we have ${(2^{i - 1})}^{n/i} = 2^{n - n/i}$ different
$\bar\sigma \in S_2^n$. Multiplying the expressions for $\tau$ and $\bar\sigma$
we get
\begin{equation}\label{eq5}
\psi_n(i) = \frac{2^n n!}{2^{n/i} \: i^{n/i} \: (n/i)!}
\end{equation}
for $i$ odd.

Let now $i \divides 2n$ be even. For $\pi = (\tau, \bar\sigma)$ to
have cycle type $i^{2n/i}$ the permutation $\tau$ must have cycle type
$(i/2)^l i^k$ with $l \ge 0, k \ge 0$, $l \cdot i/2 + k \cdot i = n$.
The number of such $\tau \in S_n$ equals
\[
\frac{n!}{(i/2)^l \, l! \: i^k \, k!}.
\]
We fix $\tau$ and count the number of corresponding $\bar\sigma$. 
Each $i/2$-cycle of $\tau$ induces one cycle of $(\tau, \bar\sigma)$ of the
length $i$. Hence $|L|$ must be odd for each $i/2$-cycle of $\tau$. There are
\[
\sum_{m \leq i/2, \; m \: \text{odd}} \binom{i/2}{m} = 2^{i/2 - 1}
\]
choices for $L$. Each $i$-cycle of $\tau$ induces two cycles of 
$(\tau, \bar\sigma)$, both having length $i$. Hence $|L|$ must be even for each
$i$-cycle of $\tau$, and we again have $2^{i - 1}$ choices for $L$. It follows, there
are ${(2^{i/2 - 1})}^l \cdot {(2^{i - 1})}^k = 2^{n - l - k}$ different
$\bar\sigma \in S_2^n$ such that $(\tau, \bar\sigma)$ has cycle type $i^{2n/i}$. 
Multiplying the expressions for $\tau$ and $\bar\sigma$ and summing up over all 
admissible $l, k$ we obtain
\begin{equation}\label{eq6}
\psi_n(i) = \sum_{\substack{l \ge 0, k \ge 0 \\ l \cdot i/2 + k \cdot i = n}} 
\frac{2^n \: n!}{2^l \: (i/2)^l \: l! \; 2^k \: i^k \: k!}
= \frac{2^n \: n!}{i^{2n/i}} \; \sum_{k = 0}^{\lfloor \frac{n}{i} \rfloor}
\frac{i^k}{(2n/i - 2k)! \: 2^k \: k!}
\end{equation}
for $i$ even.

It remains to substitute (\ref{eq5}) and (\ref{eq6}) into~(\ref{eq4}).
Setting
\[
\nu_n(i) = \frac{i^{2n/i} \: (2n/i)!}{2^n \: n!} \: \psi_n(i)
\]
we have
\[
\nu_n(i) = \frac{i^{n/i} \: (2n/i)!}{2^{n/i} \: (n/i)!} = i^{n/i} \: (2n/i - 1)!!
\]
for $i$ odd and
\[
\nu_n(i) = \sum_{k = 0}^{\lfloor \frac{n}{i} \rfloor}
\frac{i^k \: (2n/i)!}{(2n/i - 2k)! \: 2^k \: k!} = \sum_{k = 0}^{\lfloor \frac{n}{i} \rfloor}
i^k \: \binom{2n/i}{2k} \: (2k - 1)!!
\]
for $i$ even which completes the proof.
\end{proof}

As the terms in (\ref{eq2}) are all positive it is clear that 
$\underline{c}_n = (2n)^{-1} \, (2n - 1)!!$ is a lower bound for $c_n$
for $n \ge 1$. It is an easy matter to show that $\underline{c}_n$ is 
actually an asymptotic estimate for $c_n$ as $n \to \infty$.
\begin{cor}\label{c1}
\begin{equation}\label{eq7}
c_n \sim \underline{c}_n \quad \text{as $n \to \infty$}
\end{equation}
\end{cor}

\begin{proof}
Dividing out the first term in (\ref{eq2}) gives
\begin{equation}\label{eq8}
2n c_n = (2n - 1)!! + \sum_{i \divides 2n, \; i > 1} \varphi(i) \nu_n(i).
\end{equation}
We begin with establishing an upper bound for $\nu_n(i)$ using Stirling's formula
in the following form
\begin{equation*}
n! = \sqrt{2\pi n} \left( \frac{n}{e} \right)^n e^{\frac{\theta}{12n}}
\end{equation*}
where $\theta = \theta(n)$ satisfies $0 < \theta < 1$ (see \cite{Odl}). We have
\begin{equation*}
(2k - 1)!! = \sqrt{2} \left( \frac{2k}{e} \right)^k
e^{\left( \frac{1}{2}\theta_1 - \theta_2\right)\frac{1}{12k}} < 
\sqrt{2e} \left( \frac{2k}{e} \right)^k
\end{equation*}
for $k \ge 1$, $0 < \theta_1 < 1$, $0 < \theta_2 < 1$. Then we find
\begin{equation}\label{eq10}
\nu_n(i) = (2n/i - 1)!! < \sqrt{2e}\left( \frac{2n}{ei} \right)^{n/i}
\end{equation}
for $i \divides 2n$ odd. Using the following estimate
\begin{equation*}
\binom{n}{k} \le \left( \frac{ne}{k} \right)^k
\end{equation*}
(see \cite{Odl}) for the binomal coefficients we obtain
\begin{equation*}
\binom{2n/i}{2k} \, i^k \, (2k - 1)!! 
< \sqrt{2e} \left( \frac{2n^2 e}{ik} \right)^k < \sqrt{2e} \big( 2en \big)^{n/i}
\end{equation*}
for $k \in \{ 1, \ldots, \lfloor n/i \rfloor \}$ whence
\begin{equation}\label{eq11}
\nu_n(i) = 1 + \sum_{k = 1}^{\lfloor n/i \rfloor} \binom{2n/i}{2k} \, i^k \, (2k - 1)!!
< \sqrt{2e} \, n \big( 2en \big)^{n/i}
\end{equation}
for $i \divides 2n$ even. Comparing upper bounds (\ref{eq10}) and (\ref{eq11})
we conclude that
\begin{equation*}
\nu_n(i) < \overline{\nu}_n = \sqrt{2e} \, n \big( 2en \big)^{n/2}
\end{equation*}
for $i > 1$, $i \divides 2n$.
Going back to (\ref{eq8}) we see that
\begin{equation*}
\sum_{i \divides 2n, \; i > 1} \varphi(i) \nu_n(i) < 2n \, \overline{\nu}_n 
= \osmall \big( (2n - 1)!! \big)
\end{equation*}
as $n \to \infty$ and the result follows.
\end{proof}

Corollary~\ref{c1} shows that asymptotically each equivalence class contains $2n$
diagrams, which is equivalent to saying that the fraction of 1-factors of $K_{2n}$
with a nontrivial stabilizer in $C_{2n}$ tends to zero as $n \to \infty$.

\section{}

The analysis done in the proof of Theorem~\ref{t2} allows us to handle the case
of a dihedral group $D_{2n}$.
\begin{teo}
The number $d_n$ of nonequivalent $n$-diagrams under the action of a dihedral
group $D_{2n}$ equals
\begin{equation}\label{eq12}
d_n = \frac{1}{2} \big( c_n + \frac{1}{2} \big( \kappa_{n - 1} + \kappa_n \big) \big)
\end{equation}
where
\begin{equation*}
\kappa_n = \sum_{k = 0}^{\lfloor \frac{n}{2} \rfloor} \frac{n!}{k! \, (n - 2k)!}.
\end{equation*}
\end{teo}

\begin{proof}
As $C_{2n} < D_{2n}$ it follows from (\ref{eq1}) that
\begin{equation}\label{eq13}
d_n = \frac{1}{2} \left( c_n + \frac{1}{2^n \, n! \, 2n} \: \gamma_n \right)
\end{equation}
where $\gamma_n$ represents the contribution of those permutations of $D_{2n}$
which are not in $C_{2n}$. Such permutations have cycle type either $2^n$ or
$1^2 \, 2^{n - 1}$ and there are $n$ permutations of each type in $D_{2n}$.
The product in (\ref{eq1}) equals $2^n \, n!$ for permutations of cycle type $2^n$
and $2^n \, (n - 1)!$ for permutations of cycle type $1^2 \, 2^{n - 1}$. So we
can write
\begin{equation}\label{eq14}
\gamma_n = 2^n \, n! \, n \, \psi_n(2) + 2^n \, (n - 1)! \, n \, \Psi_n
\end{equation}
where $\Psi_n$ is the number of $\pi \in S_n \wr S_2$ of cycle type $1^2 \, 2^{(n - 1)}$
and $\psi_n(2)$ is the number of $\pi$ of cycle type $2^n$. Applying the analysis in the
proof of Theorem \ref{t2} we see that for each $\tau \in S_n$ of cycle type
$1^l \, 2^k$, $l + 2k = n$, $l \ge 1$, $k \ge 0$ there are $l \cdot 2^k$ permutations
$\overline\sigma \in S_2^n$ such that $\pi = (\tau, \overline\sigma) \in S_n \wr S_2^n$
has cycle type $1^2 \, 2^{(n - 1)}$. Multiplying and summing up over all admissible
$l, k$ and simplifying we get
\begin{equation*}
\Psi_n = \sum_{\substack{l \ge 1, k \ge 0 \\ l + 2k = n}} \frac{l \, n!}{l! \, k!}
= \sum_{k = 0}^{\lfloor \frac{n - 1}{2} \rfloor} \frac{n!}{(n - 1 - 2k) \, k!}.
\end{equation*}
Substituting the expressions for $\psi_n(2)$ and $\Psi_n$ into (\ref{eq14})
we obtain
\begin{equation}\label{eq15}
\gamma_n = 2^n \, n! \, n \left( 
\sum_{k = 0}^{\lfloor \frac{n}{2} \rfloor} \frac{n!}{(n - 2k)! \, k!}
+ \sum_{k = 0}^{\lfloor \frac{n - 1}{2} \rfloor} \frac{(n - 1)!}{(n - 1 - 2k)! \, k!} \right).
\end{equation}
Denoting
\begin{equation*}
\kappa_n = \sum_{k = 0}^{\lfloor \frac{n}{2} \rfloor} \frac{n!}{k! \, (n - 2k)!}
\end{equation*}
and substituting (\ref{eq15}) into (\ref{eq13}) we get the statement of the theorem.
\end{proof}

As in the case of the cyclic group, $\underline{d}_n = (4n)^{-1} \, (2n - 1)!!$ is
a lower bound for $d_n$ for $n \ge 1$ and in fact an asymptotic estimate.
\begin{cor}\label{cor2}
\begin{equation}\label{eq16}
d_n \sim \underline{d}_n \quad \text{as $n \to \infty$}
\end{equation}
\end{cor}

\begin{proof}
From (\ref{eq12}) we get
\begin{equation*}
4n d_n = 2n c_n + n (\kappa_{n - 1} + \kappa_n).
\end{equation*}
Clearly
\begin{equation*}
n (\kappa_{n - 1} + \kappa_n) < 2n\kappa_n < 
2n \, n! \sum_{k = 0}^{\lfloor \frac{n}{2} \rfloor} \frac{1}{k!} < 2n^2 \, n!
\end{equation*}
for $n \ge 1$. Stirling's formula shows that
$2n^2 \, n! = \osmall \, ((2n - 1)!!)$ and hence 
$n (\kappa_{n - 1} + \kappa_n) = \osmall \, ((2n - 1)!!)$. 
But $2n c_n \sim (2n - 1)!!$ which completes the proof.
\end{proof}

Corollary \ref{cor2} shows that asymptotically each equivalence class 
contains $4n$ diagrams.

\section{}

The following table gives an idea of the growth rate of $c_n$ and $d_n$
along with the integral parts of the corresponding asymptotic estimates.
\smallskip

\begin{center}
\begin{tabular}{|r|r|r|r|r|} \hline
\multicolumn{1}{|c|}{$n$} & 
\multicolumn{1}{c|}{$c_n$} & 
\multicolumn{1}{c|}{$\lfloor \underline{c}_n \rfloor$} & 
\multicolumn{1}{c|}{$d_n$}  & 
\multicolumn{1}{c|}{$\lfloor \underline{d}_n \rfloor$} \\ \hline
3  &  5          &  2          &  5          &  1          \\
4  &  18         &  13         &  17         &  6          \\
5  &  105        &  94         &  79         &  47         \\
6  &  902        &  866        &  554        &  433        \\
7  &  9749       &  9652       &  5283       &  4826       \\
8  &  127072     &  126689     &  65346      &  63344      \\
9  &  1915951    &  1914412    &  966156     &  957206     \\
10 &  32743182   &  32736453   &  16411700   &  16368226   \\
11 &  625002933  &  624968662  &  312702217  &  312484331  \\      \hline
\end{tabular}
\end{center}
\bigskip

Below are shown all nonequivalent (under the cyclic group) 3- and 4-diagrams.
Except for the two diagrams 12 and 13 all of them are also nonequivalent
under the dihedral group.

\bigskip

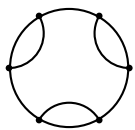
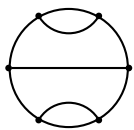
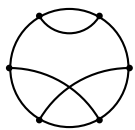
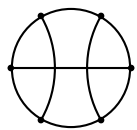
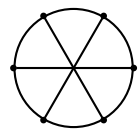
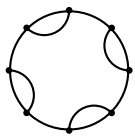
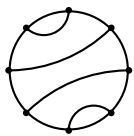
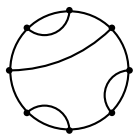
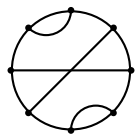
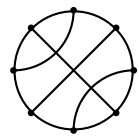
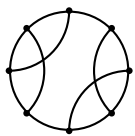
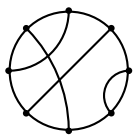
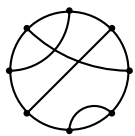
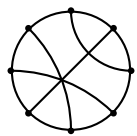
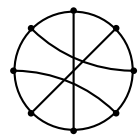
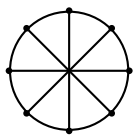
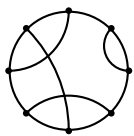
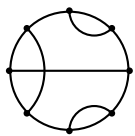
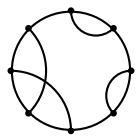
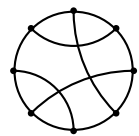
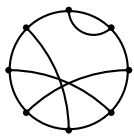
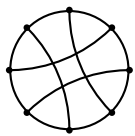
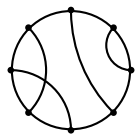
\begin{center}
\begin{pspicture}(0,0)(12.8,3.2)
\pscircle(0.800,1.600){0.80}
\psarc{*-*}(1.600,2.062){0.462}{150.0}{270.0}
\psarc{*-*}(0.000,2.062){0.462}{270.0}{30.0}
\psarc{*-*}(0.800,0.676){0.462}{30.0}{150.0}
\rput(0.8,0.2){1}
\pscircle(3.600,1.600){0.80}
\psline{*-*}(4.400,1.600)(2.800,1.600)
\psarc{*-*}(3.600,2.524){0.462}{210.0}{330.0}
\psarc{*-*}(3.600,0.676){0.462}{30.0}{150.0}
\rput(3.6,0.2){2}
\pscircle(6.400,1.600){0.80}
\psarc{*-*}(7.200,0.214){1.386}{90.0}{150.0}
\psarc{*-*}(6.400,2.524){0.462}{210.0}{330.0}
\psarc{*-*}(5.600,0.214){1.386}{30.0}{90.0}
\rput(6.4,0.2){3}
\pscircle(9.200,1.600){0.80}
\psline{*-*}(10.000,1.600)(8.400,1.600)
\psarc{*-*}(10.800,1.600){1.386}{150.0}{210.0}
\psarc{*-*}(7.600,1.600){1.386}{330.0}{30.0}
\rput(9.2,0.2){4}
\pscircle(12.000,1.600){0.80}
\psline{*-*}(12.800,1.600)(11.200,1.600)
\psline{*-*}(12.400,2.293)(11.600,0.907)
\psline{*-*}(11.600,2.293)(12.400,0.907)
\rput(12.0,0.2){5}
\end{pspicture}
\end{center}

\begin{center}
\begin{pspicture}(0,0)(12.8,3.2)
\pscircle(0.800,1.600){0.80}
\psarc{*-*}(1.600,1.931){0.331}{135.0}{270.0}
\psarc{*-*}(0.469,2.400){0.331}{225.0}{0.0}
\psarc{*-*}(0.000,1.269){0.331}{315.0}{90.0}
\psarc{*-*}(1.131,0.800){0.331}{45.0}{180.0}
\rput(0.8,0.2){6}
\pscircle(3.600,1.600){0.80}
\psarc{*-*}(4.400,-0.331){1.931}{90.0}{135.0}
\psarc{*-*}(2.800,3.531){1.931}{270.0}{315.0}
\psarc{*-*}(3.269,2.400){0.331}{225.0}{0.0}
\psarc{*-*}(3.931,0.800){0.331}{45.0}{180.0}
\rput(3.6,0.2){7}
\pscircle(6.400,1.600){0.80}
\psarc{*-*}(7.200,1.269){0.331}{90.0}{225.0}
\psarc{*-*}(5.600,3.531){1.931}{270.0}{315.0}
\psarc{*-*}(6.069,2.400){0.331}{225.0}{0.0}
\psarc{*-*}(6.069,0.800){0.331}{0.0}{135.0}
\rput(6.4,0.2){8}
\pscircle(9.200,1.600){0.80}
\psline{*-*}(10.000,1.600)(8.400,1.600)
\psline{*-*}(9.766,2.166)(8.634,1.034)
\psarc{*-*}(8.869,2.400){0.331}{225.0}{0.0}
\psarc{*-*}(9.531,0.800){0.331}{45.0}{180.0}
\rput(9.2,0.2){9}
\pscircle(12.000,1.600){0.80}
\psarc{*-*}(12.800,0.800){0.800}{90.0}{180.0}
\psline{*-*}(12.566,2.166)(11.434,1.034)
\psarc{*-*}(11.200,2.400){0.800}{270.0}{360.0}
\psline{*-*}(11.434,2.166)(12.566,1.034)
\rput(12.0,0.2){10}
\end{pspicture}
\end{center}

\begin{center}
\begin{pspicture}(0,0)(12.8,3.2)
\pscircle(0.800,1.600){0.80}
\psarc{*-*}(1.600,0.800){0.800}{90.0}{180.0}
\psarc{*-*}(1.931,1.600){0.800}{135.0}{225.0}
\psarc{*-*}(0.000,2.400){0.800}{270.0}{360.0}
\psarc{*-*}(-0.331,1.600){0.800}{315.0}{45.0}
\rput(0.8,0.2){11}
\pscircle(3.600,1.600){0.80}
\psarc{*-*}(4.400,1.269){0.331}{90.0}{225.0}
\psline{*-*}(4.166,2.166)(3.034,1.034)
\psarc{*-*}(2.800,2.400){0.800}{270.0}{360.0}
\psarc{*-*}(1.669,0.800){1.931}{0.0}{45.0}
\rput(3.6,0.2){12}
\pscircle(6.400,1.600){0.80}
\psarc{*-*}(7.200,3.531){1.931}{225.0}{270.0}
\psline{*-*}(6.966,2.166)(5.834,1.034)
\psarc{*-*}(5.600,2.400){0.800}{270.0}{360.0}
\psarc{*-*}(6.731,0.800){0.331}{45.0}{180.0}
\rput(6.4,0.2){13}
\pscircle(9.200,1.600){0.80}
\psarc{*-*}(10.000,2.400){0.800}{180.0}{270.0}
\psline{*-*}(9.766,2.166)(8.634,1.034)
\psarc{*-*}(7.269,0.800){1.931}{0.0}{45.0}
\psarc{*-*}(8.400,-0.331){1.931}{45.0}{90.0}
\rput(9.2,0.2){14}
\pscircle(12.000,1.600){0.80}
\psarc{*-*}(12.800,3.531){1.931}{225.0}{270.0}
\psline{*-*}(12.566,2.166)(11.434,1.034)
\psline{*-*}(12.000,2.400)(12.000,0.800)
\psarc{*-*}(11.200,-0.331){1.931}{45.0}{90.0}
\rput(12.0,0.2){15}
\end{pspicture}
\end{center}

\begin{center}
\begin{pspicture}(0,0)(12.8,3.2)
\pscircle(0.800,1.600){0.80}
\psline{*-*}(1.600,1.600)(0.000,1.600)
\psline{*-*}(1.366,2.166)(0.234,1.034)
\psline{*-*}(0.800,2.400)(0.800,0.800)
\psline{*-*}(0.234,2.166)(1.366,1.034)
\rput(0.8,0.2){16}
\pscircle(3.600,1.600){0.80}
\psarc{*-*}(4.400,1.931){0.331}{135.0}{270.0}
\psarc{*-*}(2.800,2.400){0.800}{270.0}{360.0}
\psarc{*-*}(1.669,0.800){1.931}{0.0}{45.0}
\psarc{*-*}(3.600,0.469){0.800}{45.0}{135.0}
\rput(3.6,0.2){17}
\pscircle(6.400,1.600){0.80}
\psline{*-*}(7.200,1.600)(5.600,1.600)
\psarc{*-*}(6.731,2.400){0.331}{180.0}{315.0}
\psarc{*-*}(5.269,1.600){0.800}{315.0}{45.0}
\psarc{*-*}(6.731,0.800){0.331}{45.0}{180.0}
\rput(6.4,0.2){18}
\pscircle(9.200,1.600){0.80}
\psarc{*-*}(10.000,1.269){0.331}{90.0}{225.0}
\psarc{*-*}(9.531,2.400){0.331}{180.0}{315.0}
\psarc{*-*}(8.069,1.600){0.800}{315.0}{45.0}
\psarc{*-*}(8.400,0.800){0.800}{0.0}{90.0}
\rput(9.2,0.2){19}
\pscircle(12.000,1.600){0.80}
\psarc{*-*}(12.800,-0.331){1.931}{90.0}{135.0}
\psarc{*-*}(12.000,2.731){0.800}{225.0}{315.0}
\psarc{*-*}(13.931,2.400){1.931}{180.0}{225.0}
\psarc{*-*}(11.200,0.800){0.800}{0.0}{90.0}
\rput(12.0,0.2){20}
\end{pspicture}
\end{center}

\begin{center}
\begin{pspicture}(0,0)(12.8,3.2)
\pscircle(3.600,1.600){0.80}
\psarc{*-*}(4.400,-0.331){1.931}{90.0}{135.0}
\psarc{*-*}(3.931,2.400){0.331}{180.0}{315.0}
\psarc{*-*}(1.669,0.800){1.931}{0.0}{45.0}
\psarc{*-*}(2.800,-0.331){1.931}{45.0}{90.0}
\rput(3.6,0.2){21}
\pscircle(6.400,1.600){0.80}
\psarc{*-*}(7.200,-0.331){1.931}{90.0}{135.0}
\psarc{*-*}(5.600,3.531){1.931}{270.0}{315.0}
\psarc{*-*}(8.331,2.400){1.931}{180.0}{225.0}
\psarc{*-*}(4.469,0.800){1.931}{0.0}{45.0}
\rput(6.4,0.2){22}
\pscircle(9.200,1.600){0.80}
\psarc{*-*}(10.000,1.931){0.331}{135.0}{270.0}
\psarc{*-*}(11.131,2.400){1.931}{180.0}{225.0}
\psarc{*-*}(8.069,1.600){0.800}{315.0}{45.0}
\psarc{*-*}(8.400,0.800){0.800}{0.0}{90.0}
\rput(9.2,0.2){23}
\end{pspicture}
\end{center}

\bigskip
\bigskip

{\bfseries Acknowledgement.} The author is grateful to Vladimir Sharko for 
posing the problem and stimulating conversations.

\bigskip

\begin{flushleft}
\scshape \small
Rofin-Sinar Laser GmbH\\
Neufeldstr. 16\\
85232 Guending\\
Germany\\
{\itshape KhruzinA@rofin-muc.de}
\end{flushleft}
 
\end{document}